\documentclass[12pt]{article}
\usepackage{amsmath,amssymb,amsthm}
\usepackage[mathscr]{eucal}
\usepackage{rotate,graphics,epsfig}
\usepackage{color}
\usepackage{hyperref}

\newtheorem{The}{Theorem}
\newtheorem{Exa}[The]{Example}

\newtheorem{Pro}[The]{Proposition}

\theoremstyle{definition}
\newtheorem{Def}[The]{Definition}
\newtheorem{Rem}[The]{Remark}
\numberwithin{equation}{section}
\numberwithin{The}{section}

\newcommand{\be}{\begin{eqnarray}}
\newcommand{\ee}{\end{eqnarray}}

\newcommand{\by}{\begin{eqnarray*}}
\newcommand{\ey}{\end{eqnarray*}}

\newcommand{\bn}{\begin{enumerate}}
\newcommand{\en}{\end{enumerate}}

\newcommand{\bi}{\begin{itemize}}
\newcommand{\ei}{\end{itemize}}

\def\frac#1#2{{#1 \over #2}}

\def\ww{\boldsymbol{w}}

\begin{document}
\title{Operator Regular Variation of Multivariate Liouville Distributions}
\author{
Haijun Li\footnote{{\small\texttt{lih@math.wsu.edu}}, Department of Mathematics and Statistics, Washington State University, 
Pullman, WA 99164, U.S.A.}
}
\date{June 2023}
\maketitle

\begin{abstract}

Operator regular variation reveals general power-law distribution tail decay phenomena using operator scaling, that includes multivariate regular variation with scalar scaling as a special case. In this paper, we show that a multivariate Liouville distribution is operator regularly varying if its driving function is univariate regularly varying. Our method focuses on operator regular variation of multivariate densities, which implies, as we also show in this paper,  operator regular variation of the multivariate distributions. This general result extends the general closure property of multivariate regular variation established by de Haan and Resnick \cite{HR1987} in 1987.


\end{abstract}

\section{Introduction}
\label{S1}

Multivariate regular variation describes power-law decay patterns for tail events, that are important in analyzing multivariate extremes
\cite{Resnick07}. A $d$-dimensional random vector $X$ is said to be multivariate regularly varying if its tail measure $\mathbb{P}(X\in tB)$, on Borel sets $B\subseteq \mathbb{R}^d$, converges vaguely to a limiting measure $\nu(B)$ with a scaling function that is univariate regularly varying. The vague convergence of multivariate regular variation can be often made stronger, for a wider class of tail events using operator scaling in place of scalar $t$. Such an extension is useful in analysis of various multi-dimensional extreme events \cite{MS01, BE07}, since it is evident, such as in finance, that vector data with heavy
tails need not have the same tail index in every direction and that it may be necessary to consider
rotated coordinate systems using operator scaling to detect variations in tail behavior \cite{MS99}. 

The goal of this paper is to establish the operator regular variation property for a multivariate Liouville distribution of $X=(X_1, \dots, X_d)$ with joint density    
\[
	f(x_1, \dots, x_d)=\kappa\, g\Big(\sum_{i=1}^dx_i\Big)\prod_{i=1}^dx_i^{a_i-1}, \ \ x_1>0, \dots, x_d>0, 
\]
where $\kappa>0$ is a constant, $a_i>0$, $i=1, \dots, d$, and $g: \mathbb{R}_+\to \mathbb{R}_+$ is known as the driving function. The theory of multivariate Liouville distributions can be found in \cite{Gupta87, Gupta91, Gupta92, Gupta95}, that also include various examples and applications. 
The multivariate regular variation of Liouville random vector $X$ seems driven by the univariate regular variation of $g(\cdot)$, and in fact, we establish in Section 3 a stronger result that $X$ is operator regularly varying under the same condition that the driving function $g(\cdot)$ is univariate regularly varying. 

Many multivariate distributions have the densities that are functions of certain norms $||\cdot||$ on $\mathbb{R}^d$ \cite{BE07}. While any norm is homogeneous, it is also often asymptotically quasi-homogeneous, and this property, together with local uniform convergence of univariate regular variation, yield the operator regular variation of multivariate distributions, such as the multivariate Liouville distribution. 
We review multivariate regular variation and operator regular variation in Section 2, and in particular, we prove that the operator regular variation of a density implies the operator regular variation of the corresponding  multivariate distribution. This result appears to be new in the literature, and is used in Section 3 to establish the operator regular variation of multivariate Liouville distributions.   

Multivariate/operator regular variation has found various applications in the multivariate extreme value analysis
\cite{Resnick07} and limiting theory \cite{MS01}, among others. Multivariate regular variation is shown to be equivalent to tail dependence of copulas \cite{LS2009, HJL12, Li2018}, that is a  fundamental property for many most useful copulas \cite{Li2009, JLN10}. Since the tail risk measures are often expressed in terms of tail densities of the multivariate copulas of underlying loss distributions \cite{JL11, SunLi2010, ZL12}, multivariate/operator regular variation has become especially useful in risk management \cite{BE07}. 
 
A univariate measurable function $f(\cdot)$ is said to be regularly varying with tail index $-\alpha$, denoted by $f\in \mbox{RV}_{-\alpha}$, if ${f(tx)}/{f(t)} \to x^{-\alpha}$, for $x>0$, where $\alpha\in \mathbb{R}$. A function $f\in \mbox{RV}_0$ is called slowly varying. 
We consider throughout this paper that any involved slow varying function is continuous. The assumption is rather mild due to 
 Karamata's representation (see, e.g., \cite{BGT1987, Resnick07}) that any slow varying function can be written as the product of a continuous function and a measurable function with  positive constant limit. All the functions, measures, and sets discussed in this paper are assumed to be measurable without explicit mention. For any two vectors in $\mathbb{R}^d$, their relations and operations, such as multivariate intervals, are taken component-wise.

\section{Distributions with operator regularly varying densities}
\label{S2}

A multivariate density function $f:\mathbb{R}^d_+\to \mathbb{R}_+$ is said to be multivariate regularly varying, denoted as $f\in \mbox{MRV}(-\rho, \lambda(x))$, if the convergence 
\begin{equation}
	\label{MRV-1}
	\frac{f(tx)}{t^{-d}V(t)}\to \lambda(x),\ \ 
\end{equation}
holds locally uniformly in $x\in \mathbb{R}^d_+\backslash\{0\}$, for $V(t)\in \mbox{RV}_{-\rho}$. Observe that $t^{-d}V(t)\in \mbox{RV}_{-d-\rho}$, and therefore the tails of a multivariate regularly varying density $f(tx)$, $x\in \mathbb{R}^d_+$, enjoy univariate power-law decays along a ray $tx$ with rate $\lambda(x)$, as $t\to \infty$. A multivariate regularly varying density of a random vector $X$ implies the multivariate regular variation of its distribution, as was shown in \cite{HR1987}.

\begin{The}\rm
	\label{tail density} (de Haan and Resnick, \cite{HR1987}) 
	Assume the density $f$ of the distribution $F$ of $X$ exists and the margins $F_i$, $1\le i\le d$, are regularly varying with tail index $\alpha>0$. If 
	$\frac{f(tx)}{t^{-d}\overline{F}_1(t)}\to \lambda(x)>0$, as $t\to \infty$,  on $\overline{\mathbb{R}}^d_+\backslash\{0\}$ and uniformly on $\{x>0: ~||x||=1\}$ where $\lambda(\cdot)$ is bounded, 
	then, for any $x\in {\mathbb{R}}^d_+\backslash\{0\}$, 
	\begin{equation}
		\label{MRV-2}
		\lim_{t\to \infty}\frac{\mathbb{P}(X\in t\,[0,x]^c)}{\overline{F}_1(t)}= 
	\lim_{t\to \infty}\frac{1-F(tx)}{\overline{F}_1(t)}=\nu([0,x]^c)=\int_{[0,x]^c}\lambda(y)dy, 
	\end{equation}
	with homogeneous property that $\lambda(tx)=t^{-\alpha-d}\lambda(x)$ for $t>0$. 
\end{The}  
The tail density of a copula is introduced in \cite{LW2013, LH14}, and using Theorem \ref{tail density}, the tail density of a copula implies the tail dependence of the copula \cite{LW2013, LH14}. 
The convergence \eqref{MRV-2} is known as the multivariate regular variation of  $X$ with limiting measure $\nu(\cdot)$. Theorem \ref{tail density} was restricted to the cone $\mathbb{R}^d_+$, but can be easily extended to the entire  $\mathbb{R}^d$ with similar proof arguments, that depend only on the algebraic structure of cones.

Theorem \ref{tail density} can be generalized if the simple univariate scaling $t$ is replaced by the operator scaling. 
Given a $d\times d$ matrix $A$, we define the exponential matrix
\[\exp(A)=\sum_{k=0}^\infty \frac{A^k}{k!}, \ \mbox{where}\ A^0=I\ \mbox{(the $d\times d$ identity matrix)}, 
\]
and the power matrix
\begin{equation}
	\label{power matrix}
	t^A=\exp(A\log t)= \sum_{k=0}^\infty \frac{A^k(\log t)^k}{k!}, \ \mbox{for}\ t>0.
\end{equation}
Power matrices can be viewed as linear operators from $\mathbb{R}^d$ to $\mathbb{R}^d$, and behave like power functions; for example, $t^{-A}=(t^{-1})^A=(t^A)^{-1}$. For any positive-definite matrix $A$ and any norm $||\cdot||$ on $\mathbb{R}^d$, $||t^A\ww||\to \infty$, as $t\to \infty$, uniformly on compact subsets of $\ww\in \mathbb{R}^d\backslash \{0\}$. A good summary on properties of exponential and power matrices, as well as operator regular variation in general can be found in \cite{MS01}. Operator regular variation for copulas has been studied in \cite{Li2018}.  

Using power matrices as scaling, operator regular variation of a multivariate density is defined as follows.

\begin{Def}
	\label{d2.1} Suppose that a non-negative random vector $(X_1, \dots, X_d)$ has a distribution $F$ on $\mathbb{R}^d$ with density $f$. The density $f$ is said to be regularly varying with operator tail index $E$, denoted as $f\in \mbox{MRV}(E, -\rho,\lambda(x))$, if the convergence
	\begin{equation}
		\label{e2.1}
		\frac{f(t^Ex)}{t^{-\text{tr}(E)}V(t)}\to \lambda(x)>0, 
	\end{equation}
holds locally uniformly in $x\in \mathbb{R}^d\backslash\{0\}$, for $V(t)\in \mbox{RV}_{-\rho}$. Here and hereafter $\text{tr}(E)$ denotes the trace of a matrix $E$. 
\end{Def}

Obviously Definition \ref{d2.1} reduces to \eqref{MRV-1} when $E=I$, the identity matrix. 
Observe that the limiting function $\lambda(\cdot)$ satisfies the scaling property that $\lambda(s^Ex) = s^{-\text{tr}(E)-\rho}\lambda(x)$, $x\in \mathbb{R}^d\backslash\{0\}$, for all $s>0$. 
Assume throughout this paper that the matrix $E$ is positive-definite, and therefore, the following spectral decomposition holds
\begin{equation}
E = O^{-1}
\begin{pmatrix}
	\lambda_1& \cdots & 0\\
	\vdots & \ddots & \vdots \\
	0 & \cdots & \lambda_d
\end{pmatrix}
O = O^{-1}DO, 
\label{diag}
\end{equation}
where $O$ is an orthogonal matrix, eigenvalues $\lambda_1, \dots, \lambda_d$ are all positive and $D=\mbox{DIAG}(\lambda_i)$ is the diagonal matrix with diagonal entries $\lambda_i$s. It follows immediately that 
\begin{equation}
	t^E=O^{-1}
	\begin{pmatrix}
		t^{\lambda_1} & \cdots & 0\\
		\vdots & \ddots & \vdots \\
		0 & \cdots & t^{\lambda_d}
	\end{pmatrix}
	O = O^{-1}t^DO,\ \ \mbox{for}\ t>0. 
	\label{dec}
\end{equation}
With these notations, a density $f(\cdot)$ is regularly varying with operator tail index $E$ and limiting function $\lambda(\cdot)$ if and only if the density $f^*(\cdot) = f(O^{-1}\cdot)$ is regularly varying with operator tail index $D$ and limiting function $\lambda(O^{-1}\cdot)$, as the following result shows. 
\begin{Pro}
	\label{positive-definite} 
	\rm Let $f$ be a density on $\mathbb{R}^d$ and $f^*(x) = f(O^{-1}x)$, $x\in \mathbb{R}^d$. Then 	
	$f\in \mbox{MRV}(E, -\rho,\lambda(x))$ if and only if $f^*\in \mbox{MRV}(D, -\rho,\lambda(O^{-1}x))$. 
\end{Pro}

\noindent
{\sl Proof.} It follows from \eqref{dec} and the fact that when $y=Ox$, 
\[
\frac{f(t^Ex)}{t^{-\text{tr}(E)}V(t)} = \frac{f(O^{-1}t^DOx)}{t^{-\text{tr}(D)}V(t)}=\frac{f^*(t^Dy)}{t^{-\text{tr}(D)}V(t)}
\]
converges to $\lambda(x)$ locally uniformly in $x$ if and only if ${f^*(t^Dy)}/(t^{-\text{tr}(D)}V(t))$ converges to $\lambda(O^{-1}y)$ locally uniformly in $y$. 
\hfill $\Box$

Without loss of generality, one can assume that $\lambda_i\ge 1$, $i=1, \dots, d$ in \eqref{e2.1}. In the case that some eigenvalues are smaller than 1,  we can use the substitution $s=t^{\lambda_{(1)}}$, where $\lambda_{(1)} = \min\{\lambda_1, \dots, \lambda_d\}>0$, then \eqref{e2.1} is equivalent to 
\[
	\frac{f(s^{E'}x)}{s^{-\text{tr}(E')}V(s^{1/\lambda_{[1]}})}\to \lambda(x)>0, \ \mbox{as}\ s\to \infty, 
\]
holds locally uniformly in $x\in \mathbb{R}^d\backslash\{0\}$, for $V(s^{1/\lambda_{(1)}})\in \mbox{RV}_{-\rho/\lambda_{(1)}}$ and a positive-definite matrix $E' = O^{-1}s^{\scriptsize\mbox{DIAG}(\lambda_i/\lambda_{(1)})}O$, having eigenvalues $\lambda_i/\lambda_{(1)}$, $i=1, \dots, d$, that are greater than or equal to 1. Similarly, one can assume that  $\lambda_i\le 1$, $i=1, \dots, d$. Observe that the scaling functions in \eqref{e2.1} are not unique and may not be tail-equivalent.

\begin{The}
	\label{pdf->CDF} \rm
	Let $X=(X_1, \dots, X_d)$ have a distribution $F$ defined on $\mathbb{R}^d$. If $F$ has a density $f\in \mbox{MRV}(E, -\rho,\lambda(t))$ on $\{x: ||x||\ge \epsilon\}$, $\epsilon>0$, for any norm $||\cdot||$ on $\mathbb{R}^d$, where $\lambda(x)$ is locally bounded and $E$ is a positive-definite matrix, then for any Borel subset $B\subseteq \{x: ||x||\ge \epsilon\}$, $\epsilon>0$,
	\begin{equation}
		\label{e2.3}
	\frac{\mathbb{P}(X\in t^EB)}{V(t)}\to \int_B\lambda(x)dx,
	\end{equation}
where $V(t)\in \mbox{RV}_{-\rho}$.
\end{The} 

\noindent
{\sl Proof.} Assume that for some $V(t)\in \mbox{RV}_{-\rho}$, the limit
\begin{equation}
	\label{e2.4}
	\frac{f(t^Ex)}{t^{-\text{tr}(E)}V(t)}\to \lambda(x)>0, 
\end{equation}
holds locally uniformly in $x\in R_\epsilon:=\{x: ||x||\ge \epsilon\}\subset \mathbb{R}^d$, $\epsilon>0$, where $||\cdot||$ is a norm on $\mathbb{R}^d$. We first prove \eqref{e2.3} when $E=\mbox{DIAG}(\lambda_i)$, where $\lambda_i>0$, $1\le i\le d$. In this case, $t^E = \mbox{DIAG}(t^{\lambda_i})$, the diagonal matrix with diagonal entries $t^{\lambda_i}$s. 

Let $[x] := \sum_{i=1}^d |x_i|^{1/\lambda_i}$, $x\in \mathbb{R}^d_+$. While $[x]$ is not a norm, the function is known as a quasi-homogeneous function, with scaling $[t^Ex] = t[x]$, $t>0$. Since $[x]$ is unbounded if one of $x_i$s goes to $\pm\infty$, the set $Q=\{x\in {\cal C}: [x]=1\}$ is compact. Therefore, the limit \eqref{e2.4} holds uniformly on $Q\cap R_\epsilon$. 

Let $h(t) = t^{-\text{tr}(E)}V(t)\in \mbox{RV}_{-\text{tr}(E)-\rho}$. Consider any Borel subset $B$ of $R_\epsilon$ that is bounded away from zero, $\epsilon>0$. Then  for $x\in B$, 
\[\frac{f(t^Ex)}{h(t)}= \frac{f(t^E[x]^E\,[x]^{-E}x)}{h(t[x])}\frac{h(t[x])}{h(t)}. 
\]
Since $[x]^{-E}x\in Q$ and \eqref{e2.4} holds uniformly on $Q\cap R_\epsilon$, for any given $\delta>0$, when $t[x]>t_1$, 
\[\frac{f\big((t[x])^E\,[x]^{-E}x\big)}{h(t[x])}\le \sup_{x\in B}\lambda\big([x]^{-E}x\big) + \delta \le \sup_{x\in Q}\lambda(x)+\delta = \kappa <\infty,
\]
which follows from the compactness of $Q$ and local boundedness of $\lambda(x)$. Let $\epsilon_0$ be the smallest value of $[x]$ on $R_\epsilon$, and obviously $\epsilon_0>0$. 
Therefore, whenever $t>t_1\epsilon_0^{-1}$, 
\begin{equation}
	\label{e2.5}
	\sup_{x\in B}\frac{f\big((t[x])^E\,[x]^{-E}x\big)}{h(t[x])}\le \sup_{x\in B}\lambda\big([x]^{-E}x\big) + \delta \le \sup_{x\in Q}\lambda(x)+\delta = \kappa <\infty.
\end{equation}
In addition, since $h(t)\in \mbox{RV}_{-\text{tr}(E)-\rho}$, we have, by Karamata's representation and the uniform convergence theorem, whenever $t>t_2$, 
\begin{equation}
	\label{e2.6}
	\frac{h(t[x])}{h(t)}\le  c[x]^{-\text{tr}(E)-\rho+\gamma}, \ x\in B,
\end{equation}	
for any small $0<\gamma<\rho$ and a constant $c>0$. Therefore, it follows from \eqref{e2.5} and \eqref{e2.6} that, whenever $t>\max\{t_1\epsilon_0^{-1},t_2\}$, $f(t^Ex)/h(t)\le \kappa_1  [x]^{-\text{tr}(E)-\rho+\gamma}$, $x\in B$, where $\kappa_1>0$ is a constant. 

To show that $[x]^{-\text{tr}(E)-\rho+\gamma}$ is Lebesque integrable on $R_\epsilon$, consider the following decomposition:
\[\{x: ||x||\ge  \epsilon\} = \mbox{I}+\mbox{II}, \  \ \mbox{I}\cap \mbox{II}=\emptyset, 
\]
where $\mbox{I} = \{x: |x_i|\le 1, 1\le i\le d\}\cap R_\epsilon$ is compact, and $\mbox{II}=\{x: |x_j|>1, j\in A; x_i\le 1, i\in A^c\}\cap R_\epsilon$, for some $\emptyset \ne A\subseteq \{1, \dots, d\}$. Obviously $[x]^{-\text{tr}(E)-\rho+\gamma}$ is bounded in $\mbox{I}$ and thus integrable on $\mbox{I}$. Observe that 
\by
& & \int_{\mbox{II}} \Big(\sum_{i=1}^d |x_i|^{1/\lambda_i}\Big)^{-\text{tr}(E)-\rho+\gamma}dx \\
&\le  & \int_{\{x: |x_j|>1, j\in A\}}\int_{\{x: 0\le |x_i|\le 1, i\in A^c\}}\Big(\sum_{i=1}^d |x_i|^{1/\lambda_i}\Big)^{-\text{tr}(E)-\rho+\gamma}dx_{A^c}dx_A\\
&\le &  2^{|A^c|}\int_{\{x: |x_j|>1, j\in A\}}\Big( \sum_{j\in A}|x_j|^{1/\lambda_j}\Big)^{-\text{tr}(E)-\rho+\gamma}dx_{A}
\ey
\by
&= & 2^{|A^c|}\int_{\{y: |y_j|>1, j\in A\}}\Big(\sum_{j\in A}|y_j|\Big)^{-\text{tr}(E)-\rho+\gamma}\prod_{j\in A}\lambda_j|y_j|^{\lambda_j-1}dy_{A}\\
&\le &  2^{|A^c|}|A|^{-\text{tr}(E)-\rho+\gamma} \int_{\{y: |y_j|>1, j\in A\}}y_{(|A|)}^{-\text{tr}(E)-\rho+\gamma} y_{(|A|)}^{\sum_{j\in A}\lambda_j-|A|}\prod_{j\in A}\lambda_j\,dy_{A}\\
&= &   2^{|A^c|}|A|^{-\text{tr}(E)-\rho+\gamma}\prod_{j\in A}\lambda_j \int_{\{y: |y_j|>1, j\in A\}}\big(y_{(|A|)}\big)^{-|A|-\text{tr}(E)+\sum_{j\in A}\lambda_j-\rho+\gamma}dy_{A}<\infty, 
\ey
which follows from the facts that $y_{(|A|)}=\max\{|y_j|, j\in A\}$ is the $L_\infty$-form on $\mathbb{R}^{|A|}$ and $-|A|-\text{tr}(E)+\sum_{j\in A}\lambda_j-\rho+\gamma<-|A|$. 

Since the function $\kappa_1 [x]^{-\text{tr}(E)-\rho+\gamma}$ is Lebesque integrable on $R_\epsilon$, it then follows from dominated convergence that $\lambda(x)$ is Lebesque integrable on $B\subseteq R_\epsilon$ and 
\[\frac{\mathbb{P}(X\in t^EB)}{V(t)}=\int_B \frac{f(t^Ex)}{t^{-\text{tr}(E)}V(t)}dx\to \int_B\lambda(x)dx
\]
and \eqref{e2.3} follows. 

In the general case of positive-definiteness, it follows from Proposition \ref{positive-definite} that the density $f(O^{-1}\cdot)$ of $OX$ is regularly varying on $O{R}_\epsilon$ and converges to $\lambda(O^{-1}\cdot)$  with operator tail index $D=\mbox{DIAG}(\lambda_i)$, where $\lambda_i$s are eigenvalues of $E$. Since any orthogonal matrix is unitary, $O{R}_\epsilon = R_\epsilon$. 
Hence,
\[\frac{\mathbb{P}(OX\in t^{D}OB)}{V(t)}\to \int_{OB}\lambda(O^{-1}y)dy
\]
for a Borel set $B\subseteq R_\epsilon$. That is,
\[\frac{\mathbb{P}(X\in t^EB)}{V(t)}=\frac{\mathbb{P}(X\in O^{-1}t^{D}OB)}{V(t)}\to \int_{OB}\lambda(O^{-1}y)dy =\int_B\lambda(x)|O|dx=\int_B\lambda(x)dx,  
\]
for any Borel set $B\subseteq R_\epsilon$. 
\hfill $\Box$

\begin{Rem}
	\label{pdf->CDF r1}
	\begin{enumerate}
		\item The univariate regular variation of $V(t)$ is used only to derive asymptotic upper bound \eqref{e2.6}. Therefore, Theorem \ref{pdf->CDF} still holds if regular variation of $V(t)$ is replaced by assuming local uniform, asymptotic upper bound in terms of certain power functions; that is, whenever $t>t_2$, 
		\begin{equation}
			\label{e2.6-1}
				\frac{h(t[x])}{h(t)}\le  c[x]^{-\text{tr}(E)-\rho+\gamma}, \ x\in B,
		\end{equation}
	for some constants $0\le \gamma<\rho$ and $c>0$. 
		\item If a density is of operator regularly variation, satisfying the local uniform condition, then, according to Theorem \ref{pdf->CDF}, the distribution is operator regularly varying in the sense of \eqref{e2.3}. The result extends the main result of \cite{HR1987} to operator regular variation, including hidden regular variation as a special case. 
		\item The intensity measure $\mu(B) : = \int_B\lambda(x)dx$ satisfies the scaling property that $\mu(t^EB) = t^{-\rho}\mu(B)$, $t>0$, for any Borel subset $B\subseteq \{x: ||x||\ge \epsilon\}$, $\epsilon>0$. 
	\end{enumerate}
\end{Rem}

\begin{Rem}
	\label{pdf->CDF r2}
	The cone $\mathbb{R}^n$ is closed under operator scaling $t^E$, where $E$ is any positive-definite $d\times d$ matrix. If $E = D$ is a diagonal matrix, then sub-cone $\mathbb{R}^n_+$ of $\mathbb{R}^n$ is closed under operator scaling $t^D$; that is,  if $B\subseteq \mathbb{R}^n_+$, then $t^DB\subseteq \mathbb{R}^n_+$. Theorem \ref{pdf->CDF} holds within $\mathbb{R}^n_+$ or any non-empty sub-cone $\cal C$ of $\mathbb{R}^n$, under operator scaling $t^D$. Some distributions allow multiple hidden regular variation properties on different cub-cones with different scalings \cite{Resnick07}, and Theorem \ref{pdf->CDF} on sub-cone $\cal C$ can be applied to these situations of hidden regular variation. 
\end{Rem}

\section{Regular variation of multivariate Liouville distributions}
\label{S3}

An absolutely continuous non-negative random vector $X=(X_1, \dots, X_d)$ is said to
have a Liouville distribution, denoted by $X\sim L_d[g(t); a_1, \dots, a_d]$,  if the its joint probability
density function is proportional to
\begin{equation}
	\label{Liouville0}
	g\Big(\sum_{i=1}^dx_i\Big)\prod_{i=1}^dx_i^{a_i-1}
\end{equation}
for $x_1>0, \dots, x_d>0$, where $a_i>0$, $i=1, \dots, d$, and the driving function $g(\cdot)$ is a suitably chosen non-negative continuous function, satisfying the integrablibity that
\begin{equation}
	\label{Liouville1}
	\int_0^\infty t^{\sum_{i=1}^da_i-1}g(t)dt < \infty. 
\end{equation}
This condition is assumed to ensure that \eqref{Liouville0} is a probability density function, due to the well-known formula for Liouville's integral. 
We also assume throughout this paper that $g(\cdot)$ has the non-compact support $[0, \infty)$. 
For example, the inverted Dirichlet distribution has the joint density function 
\[f(x_1, \dots, x_d) = \frac{\Gamma\big(\sum_{i=1}^{d+1}a_i\big)}{\Gamma(a_{d+1})}\Big(1+\sum_{i=1}^dx_i\Big)^{-a_1- \cdots -a_d-a_{d+1}}\prod_{i=1}^d\frac{x_i^{a_i-1}}{\Gamma(a_i)}, 
\]
in which, $g(t) = (1+t)^{-a_1-\dots -a_d-a_{d+1}}$, $t>0$, $a_i>0$, $i=1, \dots, d+1$, where $\Gamma(\cdot)$ denotes the gamma function. Observe that this function is regularly varying with tail index $\sum_{i=1}^{d+1}a_i$. In general, however, the function $f(\cdot)$ can be any non-negative function, including rapidly varying functions.  

In general, the process of conditioning a random vector on the sum of its components leading a distribution of certain independent events can be modeled using the multivariate Liouville distribution and its various extensions. The theory and extensive discussions of multivariate Liouville distributions can be found in \cite{Gupta87, Gupta91, Gupta92, Gupta95}, and the history and related references are detailed in \cite{Gupta01}. The multivariate Liouville distributions can be extended to locally compact Abelian groups \cite{Gupta95}, that include the space of real symmetric positive-definite matrices \cite{Gupta87}. 

The univariate regular variation of the driving function $g(\cdot)$ naturally implies the joint multivariate regular variation of a multivariate Liouville distribution, as the following result shows.

\begin{Pro}
	\label{Liouville2}\rm
	If $g\in \mbox{RV}_{-\beta}$, then $X$ is regularly varying with limiting measure
	\[\mu(B) =\int_B \Big(\sum_{i=1}^dx_i\Big)^{-\beta}\prod_{i=1}^dx_i^{a_i-1}dx
	\]
	for any Borel subset $B\subseteq \mathbb{R}^d_+$ that is bounded away from 0.  
\end{Pro}

\noindent
{\sl Proof.} Since $g\in \mbox{RV}_{-\beta}$, we have
\[\frac{g\Big(\sum_{i=1}^dtx_i\Big)}{g(t)}\to \Big(\sum_{i=1}^dx_i\Big)^{-\beta},
\]
locally uniformly. 
Let $V(t) = g(t)t^{\sum_{i=1}^da_i}$, and obviously $V\in \mbox{RV}_{-\beta+\sum_{i=1}^da_i}$, implying that
\[\frac{f(tx_1, \dots, tx_d)}{t^{-d}V(t)}\to  \Big(\sum_{i=1}^dx_i\Big)^{-\beta}\prod_{i=1}^dx_i^{a_i-1},
\]
locally uniformly. Note that the condition $g\in \mbox{RV}_{-\beta}$ with \eqref{Liouville1} actually implies that $\beta>\sum_{i=1}^da_i$, and thus $V(t)\to 0$ as $t\to \infty$. 
It then follows from Theorem \ref{pdf->CDF} that 
\[	\frac{\mathbb{P}(X\in tB)}{V(t)}\to \int_B \Big(\sum_{i=1}^dx_i\Big)^{-\beta}\prod_{i=1}^dx_i^{a_i-1}dx
\]
for any Borel subset $B\subseteq \mathbb{R}^d_+$ that is bounded away from 0. 
\hfill $\Box$

A stronger result in fact holds under the same condition that $g(\cdot)$ is regularly varying.

\begin{The}
	\label{Liouville3}\rm
	If $g\in \mbox{RV}_{-\beta}$, then $X$ is regularly varying with operator tail index $E=\mbox{DIAG}\,(\alpha_i)$, for any $\alpha_i>0$, $i=1, \dots, d$, and with limiting measure
	\[\mu(B) =\int_B\Big(\sum_{i\in (\alpha)} x_i\Big)^{-\beta}\prod_{i=1}^dx_i^{a_i-1}dx
	\]
	for any Borel subset $B\subseteq \mathbb{R}^d_+$ that is bounded away from 0, where $(\alpha) = \{i: \alpha_i = \max_{1\le k\le d}\{\alpha_k\}\}$.  
\end{The}

\noindent
{\sl Proof.} Let $\alpha=\max_{1\le k\le d}\{\alpha_k\}$ and $(\alpha)=\{i: a_i=\alpha\}$, where $\alpha_i>0$, $i=1, \dots, d$. 

Observe that $\sum_{i=1}^dt^{\alpha_i-\alpha}x_i \to \sum_{i\in (\alpha)} x_i$, as $t\to \infty$. Since $g\in \mbox{RV}_{-\beta}$, it follows from the local uniform convergence that 
\[\frac{g\Big(\sum_{i=1}^dt^{\alpha_i}x_i\Big)}{g(t^\alpha)}=\frac{g\Big(t^\alpha\sum_{i=1}^dt^{\alpha_i-\alpha}x_i\Big)}{g(t^\alpha)}\to \Big(\sum_{i\in (\alpha)} x_i\Big)^{-\beta},
\]
locally uniformly. 
Let $V(t) = g(t^\alpha)t^{\sum_{i=1}^d\alpha_ia_i}$, and obviously $V\in \mbox{RV}_{-\alpha\beta+\sum_{i=1}^d\alpha_ia_i}$. It follows from \eqref{Liouville1} that $\beta>\sum_{i=1}^da_i$, implying that $\alpha\beta>\sum_{i=1}^d\alpha a_i\ge \sum_{i=1}^d\alpha_ia_i$. Therefore $V(t)\to 0$ as $t\to \infty$. Since 
\[\frac{f(t^{\alpha_1}x_1, \dots, t^{\alpha_d}x_d)}{t^{-\sum_{i=1}^d\alpha_i}V(t)}\to  \Big(\sum_{i\in (\alpha)} x_i\Big)^{-\beta}\prod_{i=1}^dx_i^{a_i-1},
\]
locally uniformly. 
It then follows from Theorem \ref{pdf->CDF} that 
\[	\frac{\mathbb{P}(X\in t^EB)}{V(t)}\to \int_B \Big(\sum_{i\in (\alpha)} x_i\Big)^{-\beta}\prod_{i=1}^dx_i^{a_i-1}dx =: \mu(B)
\]
for any Borel subset $B\subseteq \mathbb{R}^d_+$ that is bounded away from 0.
\hfill $\Box$

\begin{Rem}
	\begin{enumerate}
		\item When $\alpha_i = 1$ for all $i=1, \dots, d$, $(\alpha)=\{1, \dots, d\}$ and Theorem \ref{Liouville3} reduces to Proposition \ref{Liouville2}. 
		\item The intensity measure satisfies the scaling property that $\mu(t^EB) = t^{-\alpha\beta+\sum_{i=1}^d\alpha_ia_i}\mu(B)$, $t>0$, for any Borel subset $B\subseteq \mathbb{R}^d_+$ that is bounded away from 0. 
		\item Theorem \ref{Liouville3} holds for a more general case of the Liouville distribution, where the joint probability
		density function is proportional to
		\begin{equation}
			\label{Liouville-1}
			g\Big(\sum_{i=1}^dx_i\Big)\prod_{i=1}^d\mu_{a_i}(x_i), \ \ x_1>0, \dots, x_d>0, 
		\end{equation}
		where measures $\mu_{a_i}(\cdot)$ is univariate regularly varying with tail index $a_i$ (see, e.g., \cite{Resnick07} for details on regularly varying measures),  $a_i>0$, $i=1, \dots, d$, and $g(\cdot)\in \mbox{RV}_{-\beta}$ is a suitably chosen non-negative continuous function with \eqref{Liouville1}. The distributions  \eqref{Liouville-1} (see \cite{Gupta95}) reinforce the idea of multivariate Liouville distributions that conditioning a random vector on the sum of its components leads a distribution of certain independent univariate distributions. 
	\end{enumerate}
\end{Rem}

Let $g: \mathbb{R}_+ \to \mathbb{R}$ be Borel-measurable. Define the Weyl fractional integral of order $\alpha>0$ as follows,
\begin{equation}
	\label{Weyl}
	W^\alpha g(t) = \frac{1}{\Gamma (\alpha)}\int_t^\infty (s-t)^{\alpha - 1}g(s)ds,\ \ \ t>0,
\end{equation}	
whenever the integral exists, where $\Gamma(\alpha)$ is the gamma function of $\alpha$. The Weyl fractional integral has been used in analyzing conditional distributions of a multivariate Liouville distribution \cite{Gupta87}. 
The following result extends  Karamata's theorem to Weyl fractional integrals. 
\begin{The}
	\label{Weyl0}\rm 
	If $g\in \mbox{RV}_{-\beta}$, then $W^\alpha g(t)\in \mbox{RV}_{\alpha-\beta}$, $\alpha<\beta$. 
\end{The}	

\noindent
{\sl Proof.} If $\alpha = 1$, then the result follows immediately from Karamata's theorem. In general, the substitution via $s=xt$ leads to 
\[W^\alpha g(t) = \frac{1}{\Gamma (\alpha)}\int_t^\infty (s-t)^{\alpha - 1}g(s)ds = \frac{t^\alpha}{\Gamma (\alpha)} \int_1^\infty (x-1)^{\alpha-1}g(xt)dx. 
\]
Consider
\begin{equation}
	\label{Weyl1}
\frac{W^\alpha g(t)}{t^\alpha g(t)} = \frac{1}{\Gamma (\alpha)}\int_1^\infty (x-1)^{\alpha-1}\frac{g(xt)}{g(t)}dx. 
\end{equation}
We now show that the ratio in \eqref{Weyl1} has a constant limit as $t\to \infty$. Since  $g\in \mbox{RV}_{-\beta}$, the uniform convergence theorem yields that for any small $\epsilon >0$, there exists a constant $N$ such that whenever $t>N$
\[(1-\epsilon)x^{-\beta} \le \frac{g(xt)}{g(t)}\le (1+\epsilon) x^{-\beta},\ \ x\in [1,\infty),
\]
which implies that 
\[(1-\epsilon)\int_1^\infty (x-1)^{\alpha-1}x^{-\beta}dx \le \int_1^\infty (x-1)^{\alpha-1}\frac{g(xt)}{g(t)}dx\le (1+\epsilon)\int_1^\infty (x-1)^{\alpha-1} x^{-\beta}dx.
\]
Let $A = \int_1^\infty (x-1)^{\alpha-1}x^{-\beta}dx$. Then integration by parts and $\alpha < \beta$ lead to
	\[A = \frac{1}{\alpha}\int_1^\infty x^{-\beta}d(x-1)^\alpha = \frac{\beta}{\alpha}\int_1^\infty (x-1)^\alpha x^{-\beta-1}dx < \frac{\beta}{\alpha} \int_1^\infty x^{\alpha-\beta -1}dx < \infty. 
	\]
Since $A$ is finite, we then have, for any $\epsilon>0$, 
\[(1-\epsilon)A \le \underline{\lim}_{t \to \infty}\int_1^\infty (x-1)^{\alpha-1}\frac{g(xt)}{g(t)}dx
\le \overline{\lim}_{t \to \infty}\int_1^\infty (x-1)^{\alpha-1}\frac{g(xt)}{g(t)}dx \le (1+\epsilon)A. 
\]
Letting $\epsilon\to 0$ leads to $\lim_{t \to \infty}\int_1^\infty (x-1)^{\alpha-1}\frac{g(xt)}{g(t)}dx = A$ exists, with a finite constant limit. It then follows from \eqref{Weyl1} that $\lim_{t \to \infty}{W^\alpha g(t)}/(t^\alpha g(t)) = A/\Gamma(\alpha) <\infty$. Since $g\in \mbox{RV}_{-\beta}$, $W^\alpha g(t)\in \mbox{RV}_{\alpha-\beta}$. 
\hfill $\Box$

It follows from \cite{Gupta87} that if $X = (X_1, \dots, X_d)\sim L_d[g(t); a_1, \dots, a_d]$ then the multivariate margin $(X_1, \dots, X_r)\sim L_r[W^ag(t); a_1, \dots, a_r]$, where $a=\sum_{i=r+1}^da_i$, $r<d$. According to Theorems \ref{pdf->CDF} and \ref{Weyl0}, if $g\in \mbox{RV}_{-\beta}$, then   $W^a g(t)\in \mbox{RV}_{a-\beta}$, implying that the multivariate margin $(X_1, \dots, X_r)$ is regularly varying, which is consistent with the fact that $(X_1, \dots, X_d)$ is regularly varying in this case. Furthermore, conditional distributions of a multivariate Liouville distribution is also regularly varying. 

\begin{The}\rm 
	\label{cond Liouville}
	Suppose that $X=(X_1, \dots, X_d)\sim L_d[g(t); a_1, \dots, a_d]$ with $g \in \mbox{RV}_{-\beta}$. Let
	\[V(B, t) = \mathbb{P}\big((X_{r+1}, \dots, X_d)\in B\ \big|\,\{X_1=tx_1, \dots, X_r=tx_r\}\big)
	\]
	for any fixed $(x_1, \dots, x_r)$, where $B\subseteq \mathbb{R}^{d-r}_+$ that is bounded away from 0.
	\begin{enumerate}
		\item The conditional distribution of $(X_{r+1}, \dots, X_d)$ given $\{X_1=x_1, \dots, X_r=x_r\}$ is regularly varying, for any fixed $(x_1, \dots, x_r)$. 
		\item There exists a function $V\in \mbox{RV}_{-\beta+\sum_{i=r+1}^dx_i}$, such that $V(tB,t)/V(t)$ converges vaguely to 
		\[	\nu_r(B)=\int_B\kappa \left(\frac{\sum_{i=1}^dx_i}{1+\sum_{i=1}^rx_i}\right)^{-\beta}\prod_{i=1}^dx_i^{a_i-1}dx
		\]
		for any Borel subset $B\subseteq \mathbb{R}^{d-r}_+$ that is bounded away from 0.
	\end{enumerate}	
\end{The} 

\noindent
{\sl Proof.} It follows from Corollary 4.3 of \cite{Gupta87} that 
\begin{equation}
	\label{conditioning}
	\big[(X_{r+1}, \dots, X_d)\big|X_1=x_1, \dots, X_r=x_r\big] \sim L_{d-r}[g_r(t); a_{r+1}, \dots, a_d]
\end{equation}
where $g_r(t) = g\big(t+\sum_{i=1}^rx_i\big)/W^ag\big(\sum_{i=1}^rx_i\big)$, $a=\sum_{i=r+1}^da_i$. 

(1) For fixed $(x_1, \dots, x_r)$, $g_r(t)\in \mbox{RV}_{-\beta}$. It then follows from Proposition \ref{Liouville2} that the conditional distribution \eqref{conditioning} is regularly varying with intensity measure
\[\mu_r(B) =\int_B \Big(\sum_{i=r+1}^dx_i\Big)^{-\beta}\prod_{i=r+1}^dx_i^{a_i-1}dx
\]
for any Borel subset $B\subseteq \mathbb{R}^{d-r}_+$ that is bounded away from 0. 

(2) The conditional density of $(X_{r+1}, \dots, X_d)$ given $\{X_1=tx_1, \dots, X_r=tx_r\}$ is proportional to 
\[g_{r,t}\Big(\sum_{i=r+1}^dx_i\Big)\prod_{i=r+1}^dx_i^{a_i-1}, 
\]
where $g_{r,t}(s) = f\big(s+\sum_{i=1}^rtx_i\big)/W^af\big(\sum_{i=1}^rtx_i\big)$, $a=\sum_{i=r+1}^da_i$. 
Observe that
\begin{equation}
	\label{g function}
\frac{g_{r,t}\Big(t\sum_{i=r+1}^dx_i\Big)}{g_{r,t}(t)}=\frac{g\Big(t\sum_{i=1}^dx_i\Big)}{g\Big(t(1+\sum_{i=1}^rx_i)\Big)} \to \left(\frac{\sum_{i=1}^dx_i}{1+\sum_{i=1}^rx_i}\right)^{-\beta}< \frac{\big(\sum_{i=r+1}^dx_i\big)^{-\beta}}{\big(2+\sum_{i=1}^rx_i\big)^{-\beta}},
\end{equation}
locally uniformly.  
Let $h(s) = s^{-(d-r)}g_{r,t}(s)s^{\sum_{i=r+1}^da_i}$, and it then follows that $h(tx)/h(t) < cx^{-\beta+\sum_{i=r+1}^da_i-d+r}$, for some constant $c= \big(2+\sum_{i=1}^rx_i\big)^{\beta}>0$, as $t\to \infty$, locally uniformly.  
The conditional density $f_r(\cdot)$ of $(X_{r+1}, \dots, X_d)$, given $\{X_1=tx_1, \dots, X_r=tx_r\}$, with scaling function $t$,  satisfies 
\[\frac{f_r(tx_{r+1}, \dots, tx_d)}{h(t)}\to \kappa \left(\frac{\sum_{i=1}^dx_i}{1+\sum_{i=1}^rx_i}\right)^{-\beta}\prod_{i=r+1}^dx_i^{a_i-1},\ \mbox{as}\ t\to \infty, 
\]
locally uniformly, where $\kappa>0$ is a constant. Note that the condition $g\in \mbox{RV}_{-\beta}$ with \eqref{Liouville1} actually implies that $\beta>\sum_{i=1}^da_i$, and thus $V(t)=g_{r,t}(t)t^{\sum_{i=r+1}^da_i}\to 0$ as $t\to \infty$. 
It then follows from Theorem \ref{pdf->CDF} and Remark \ref{pdf->CDF r1} (1) that 
\begin{equation}
	\label{cond V}
	\frac{\mathbb{P}\big((X_{r+1}, \dots, X_d)\in tB\,\big|\,X_1=tx_1, \dots, X_r=tx_r\big)}{V(t)}\to \int_B\kappa \left(\frac{\sum_{i=1}^dx_i}{1+\sum_{i=1}^rx_i}\right)^{-\beta}\prod_{i=1}^dx_i^{a_i-1}dx
\end{equation}
for any Borel subset $B\subseteq \mathbb{R}^{d-r}_+$ that is bounded away from 0, and (2) follows. 
\hfill $\Box$

\begin{Rem}
	\begin{enumerate}
		\item Theorem \ref{cond Liouville} (1) shows a multivariate margin conditioning on fixed values for its complement is multivariate regular varying. In contrast, if conditioning variables are allowed to be variable, then the scaling variable $t$ can be chosen for both conditioning and conditioned variables, so that the conditional probability $V(tB, t)$ is regularly varying, as illustrated in \eqref{cond V}.  
		\item 	 Note that the function $g\big(s+\sum_{i=1}^rx_i\big)/W^ag\big(\sum_{i=1}^rx_i\big)$, $a=\sum_{i=r+1}^da_i$, is a multivariate regularly varying function of $(s,x_1, \dots, x_r)$, as is showed in \eqref{g function}, implying that  $g_{r,t}(s)$ satisfies some regular variation property only if the scaling variable is the same as $t$. Theorem \ref{pdf->CDF} can still be applied due to the asymptotic upper bound \eqref{e2.6-1} (see Remark \ref{pdf->CDF r1} (1)).  
	\end{enumerate}
\end{Rem}

The following result on conditional distributions involves non-trivial use of Theorem \ref{Weyl0}.

\begin{The}\rm 
	\label{Liouville sum}
	Suppose that $X=(X_1, \dots, X_d)\sim L_d[g(t); a_1, \dots, a_d]$ with $g \in \mbox{RV}_{-\beta}$, $1\le r<d$, and the expectations exist.  
	\begin{enumerate}
		\item The conditional joint moment \[E\Big(\prod_{i=r+1}^dX_i^{j_i}\Big| \sum_{i=1}^rX_i=t\Big)\in \mbox{RV}_{j},
		\]
		where $j=\sum_{i=r+1}^dj_i$. 
		\item For any $h(t)\in \mbox{RV}_{-\gamma}$, for which the relevant expectations
		exist,
		\[E\Big(h\big(\sum_{i=r+1}^dX_i\big)\Big| \sum_{i=1}^rX_i=t\Big)\in  \mbox{RV}_{-\gamma},
		\]
		where $a=\sum_{i=r+1}^da_i$. 
	\end{enumerate}	
\end{The} 

\noindent
{\sl Proof.} (1) It follows from \cite{Gupta87} that
\[E\Big(\prod_{i=r+1}^dX_i^{j_i}\Big| \sum_{i=1}^rX_i=t\Big) = c\frac{W^{j+a}g(t)}{W^ag(t)}
\]
for some constant $c>0$, where $a=\sum_{i=r+1}^da_i$ and $j=\sum_{i=r+1}^dj_i$. Theorem \eqref{Weyl0} implies that $W^{j+a}g(t)\in \mbox{RV}_{j+a-\beta}$ and $W^ag(t)\in \mbox{RV}_{a-\beta}$, and the result follows. 

(2) It follows from \cite{Gupta87} that
\[E\Big(h\big(\sum_{i=r+1}^dX_i\big)\Big| \sum_{i=1}^rX_i=t\Big)\,W^ag(t) = c\int_t^\infty (y-t)^{a-1}h(y-t)g(y)dy
\]
for some constant $c$, where $a=\sum_{i=r+1}^da_i$. Consider the integral on the right-hand side, with the substitution via $y=xt$:
\[\int_t^\infty (y-t)^{a-1}h(y-t)g(y)dy = t^a\int_1^\infty (x-1)^{a-1}h(t(x-1))g(tx)dx
\]
\[= t^ah(t)g(t)\int_1^\infty (x-1)^{a-1}\frac{h(t(x-1))g(tx)}{h(t)g(t)}dx.
\]
Since $g \in \mbox{RV}_{-\beta}$ and $h(t)\in \mbox{RV}_{-\gamma}$, 
\[\frac{h(t(x-1))g(tx)}{h(t)g(t)}\to (x-1)^{-\gamma}x^{-\beta}
\]
uniformly on $[1,\infty)$, using the similar proof as the limit for the integral, as $t\to \infty$, in \eqref{Weyl1}, yields that, as $t\to \infty$, 
\[\int_1^\infty (x-1)^{a-1}\frac{h(t(x-1))g(tx)}{h(t)g(t)}dx\to c' < \infty,
\]
where $c'$ is a constant. Therefore,
\[\frac{E\Big(h\big(\sum_{i=r+1}^dX_i\big)\Big| \sum_{i=1}^rX_i=t\Big)\,W^ag(t)}{t^ah(t)g(t)}\to cc',\ \ t\to \infty. 
\]
Since $W^ag(t)\in \mbox{RV}_{\alpha-\beta}$, we have that $E\big(h\big(\sum_{i=r+1}^dX_i\big)\big| \sum_{i=1}^rX_i=t\big)\in \mbox{RV}_{-\gamma}$. 
\hfill $\Box$

Thereon \ref{Liouville sum} (2) reveals an intriguing tail phenomenon for random vectors $(X_1, \dots, X_d)$ with multivariate Liouville distributions. Since $(X_1, \dots, X_d)$ is conditional independent given the sum $\sum_{i=1}^dX_i$, $(X_1, \dots, X_d)$ is positively associated, implying that $(\sum_{i=1}^rX_i,\sum_{i=r+1}^dX_i)$ is positively associated, but whether or not $(\sum_{i=1}^rX_i,\sum_{i=r+1}^dX_i)$ is stochastically increasing is unknown. Thereon \ref{Liouville sum} (2) shows that $(\sum_{i=1}^rX_i,\sum_{i=r+1}^dX_i)$ is stochastically regularly varying in the following sense
	\[E\Big(h\big(\sum_{i=r+1}^dX_i\big)\Big| \sum_{i=1}^rX_i=t\Big)\in  \mbox{RV}_{-\gamma},
\]
whenever $h(t)$ is regularly varying.

\section{Concluding remarks}
\label{S4}

It is important to have good criteria in terms of densities which imply
the regular variation of distribution tails \cite{HR1987}, since many useful multivariate distributions lack explicit expressions and are specified only by their density functions. We obtain in this paper the condition on densities that implies the operator regular variation of the underlying multivariate distributions. We then apply this result to multivariate Liouville distributions with the densities that enjoy the natural property of asymptomatic quasi-homogeneity, provided that the driving function is univariate regularly varying. 

The vague convergence for a regularly varying multivariate distribution often holds with several, different limiting measures that enjoy a wide variety of  scaling properties. The multivariate Liouville distribution is one of such distributions and clearly shows such a diverse set of multivariate regular variations, that are useful in the analysis of various multivariate extremes.


\end{document}